\documentclass[12pt]{article}
\usepackage{amsmath,amssymb}
\usepackage{verbatim}
\usepackage{dsfont,bm}
\usepackage{color}
\usepackage{graphicx}
\usepackage{amsthm}
\usepackage{enumerate}
\usepackage{esint} %\fint

\usepackage{xcolor}
\definecolor{darkgreen}{rgb}{0,0.75,0}
\definecolor{darkred}{rgb}{0.75,0,0}
\definecolor{darkmagenta}{rgb}{0.5,0,0.5}
\usepackage[colorlinks,citecolor=darkgreen,linkcolor=darkred,urlcolor=darkmagenta]{hyperref}
\newtheorem{theorem}{Theorem}[section]
\newtheorem{thm}[theorem]{Theorem}

\newtheorem{cor}[theorem]{Corollary}
\newtheorem{lemma}[theorem]{Lemma}
\newtheorem{lem}[theorem]{Lemma}

\newtheorem{definition}[theorem]{Definition}

\newtheorem{remark}[theorem]{Remark}

\numberwithin{equation}{section}

\def\be{\begin{equation}}
\def\ee{\end{equation}}
\def\bes{\begin{equation*}}
\def\ees{\end{equation*}}

\newcommand{\mr}[1]{{\tt \href{http://www.ams.org/mathscinet-getitem?mr=#1}{MR#1}}}

\newcommand{\set}[1]{\left\{ #1 \right\}}

\newcommand{\abs}[1]{{\left\vert\kern-0.25ex #1
		\kern-0.25ex\right\vert}}
\newcommand\norm[1]{\left\lVert#1\right\rVert} %norm
\newcommand{\one}{\mathds{1}} %indicator
\newcommand{\loc}[0]{\operatorname{loc}}

  \def\sC {{\mathcal C}}
 \def\sE {{\mathcal E}} \def\sF {{\mathcal F}}
  
  \def\sL {{\mathcal L}}

 \def\bN {{\mathbb N}} 
  \def\bR {{\mathbb R}}

\def\ignore#1{}

\def\ol{\overline}           

\def\lam {\lambda}  
\def\eps{\varepsilon}

\def\Gam{\Gamma} \def\gam{\gamma}

\def\to {\rightarrow}

\def\q{\quad} \def\qq{\qquad}
\def\dint{\int\kern-.6em\int}

\newcommand\restr[2]{{% we make the whole thing an ordinary symbol
		\left.\kern-\nulldelimiterspace % automatically resize the bar with \right
		#1 % the function
		\vphantom{\big|} % pretend it's a little taller at normal size
		\right|_{#2} % this is the delimiter
	}} %restriction of a function
	\def\diam{{\mathop{{\rm diam }}}}

	\def\supp{\mathop{{\rm supp}}}
	
	\def\Cap{\operatorname{Cap}}

	\newcommand{\on}[1]{\operatorname{ #1}}
	
	\def\be{\begin{equation}}
	\def\ee{\end{equation}}
	\def\bes{\begin{equation*}}
	\def\ees{\end{equation*}}
	\def\ba{\begin{align}}
	\def\ea{\end{align}}
	\def\xxea{\end{align}}
\def\bas{\begin{align*}}
\def\eas{\end{align*}}

\def\proof{{\smallskip\noindent {\em Proof. }}}
\def\qed{{\hfill $\square$ \bigskip}}

\definecolor{dgreen}{rgb}{0, 0.6, 0.1}
\definecolor{dblue}{rgb}{0, 0.0, 0.6}
\definecolor{vdblue}{rgb}{0,.08, 0.45}
\definecolor{dred}{rgb}{0.7, 0.0, 0.0}
\definecolor{vdblue}{rgb}{0,.08, 0.45}

\definecolor{purple}{rgb}{0.6, 0.0, 0.6}
\definecolor{mytext}{rgb}{0.1, 0.1, 0.1}

\begin{document}
	
	\font\titlefont=cmbx14 scaled\magstep1
	\title{\titlefont On the length of chains in a metric space.}    
	\author{	Mathav Murugan\thanks{Research partially supported by NSERC (Canada).} 
		}
	\maketitle
	\vspace{-0.5cm}
%\ver
\begin{abstract}
We obtain an upper bound on the minimal number of points in an $\eps$-chain joining two points in a metric space.
This generalizes a bound due to Hambly and Kumagai (1999) for the case of resistance metric on certain self-similar fractals.
As an application, we deduce a condition on $\eps$-chains introduced by Grigor'yan and Telcs (2012). 
This allows us to obtain sharp bounds on the heat kernel for spaces satisfying the parabolic Harnack inequality without assuming further conditions on the metric. 
A snowflake transform on the Euclidean space shows that our bound is sharp.

\vskip.2cm
%\noindent {\it Keywords:} 
\
\end{abstract}

\section{Introduction}

The fundamental solution of the heat equation (or \emph{heat kernel}) on $\bR^n$ is given by the  Gauss Weierstrass kernel
\[
p_t(x,y)= \frac{1}{(4 \pi t)^{n/2}} \exp \left(- \frac{d(x,y)^2}{4t}\right), \q \mbox{for all $x,y \in \bR^n, t>0$.}
\]
From a probabilistic viewpoint, the heat kernel can be interpreted as the transition probability density of the diffusion generated by the Laplacian $\Delta$. More generally, for any uniformly elliptic, divergence form operator $\sL u= \sum_{i,j=0}^{n} \frac{\partial }{\partial x_i} \left(a_{ij}(x)  \frac{\partial u}{\partial x_j}\right)$ on $\bR^n$, Aronson \cite{Ar67} proved that the heat kernel $p_t(x,y)$ of the corresponding heat equation $\partial_t u - \sL u =0$ satisfies
\[
\frac{c_1}{V(x,t^{1/2})} \exp\left(- \frac{d(x,y)^2}{c_1t}\right) \le p_t(x,y) \le 
\frac{C_1}{V(x,t^{1/2})} \exp\left(- \frac{d(x,y)^2}{C_1t}\right),
\]
for all  $x,y \in \bR^n, t>0$, where $c_1,C_1 \in (0,\infty)$.
Here $V(x,r)$ denotes the Lebesgue measure of the Euclidean ball $B(x,r)$ centered at $x$ with radius $r$.

To prove the above lower bound on $p_t(x,y)$, the first step is to obtain the lower bound under the additional restriction that $d(x,y) \le C t^{1/2}$. Such an estimate is called a \emph{near diagonal lower bound}. In order to obtain the full lower bound from a near diagonal lower bound, one chooses a sequence of points (called a \emph{chain}) $\set{x_i}_{i=0}^n$ such that $x_0=x, x_n=y$ and $n \in \bN$ such that $d(x_i,x_{i+1}) \le C (t/n)^{1/2}/2$ for all $i=0,1,\ldots,n-1$. Then we use Chapman-Kolmogorov equation to obtain the estimate
\[
p_t(x,y)\ge  \int_{B(x_{n-1},C (t/n)^{1/2}/2)} \ldots \int_{B(x_1,,C (t/n)^{1/2}/2)} \Pi_{i=0}^{n-1} p_{t/n}(y_i,y_{i+1}) \,dy_1\ldots dy_{n-1},
\]
where $y_0=x,y_n=y$. By optimizing over $n$ and the sequence $\set{x_i}_{i=0}^n$ and using the near diagonal lower bound, we obtain the full lower bound on the heat kernel $p_t(x,y)$. This method of obtaining full heat kernel lower bound is called the \emph{chaining argument}.

The use of chaining argument to obtain heat kernel estimates is classical \cite{AS67,Ar67}.
Such chaining arguments are also used to obtain heat kernel estimates on fractals; see \cite{Bar98} for a general introduction to diffusions on fractals.
In this work, we address the natural converse question: {\em Do heat kernel estimates imply the existence of short chains?} Our main result provides an upper bound  on the length of chains, which in some sense is the best possible. The goal of this work is to obtain sharp quantitative bounds on the connectivity of a metric space. This will be expressed as bounds on the length of chains.
\subsection{Framework and definitions}
We recall the definition of a chain in a metric space $(X,d)$.
\begin{definition}{\rm 
 We say that a sequence  $\set{x_i}_{i=0}^N$  of points in $X$ is an \emph{$\eps$-chain} between points $x,y \in X$ if 
 \[
 x_0=x, \q x_N=y, \q \mbox{ and } \q d(x_i,x_{i+1}) < \eps \q \mbox{ for all $i=0,1,\ldots,N-1$.}
 \]
 For any $\eps>0$ and $x,y \in X$, define
 \[
 d_\eps(x,y) = \inf_{\set{x_i} \mbox{ is $\eps$-chain}} \sum_{i=0}^{N-1} d(x_i,x_{i+1}),
 \]
 where the infimum is taken over all $\eps$-chains  $\set{x_i}_{i=0}^N$ between $x,y$ with arbitrary $N$.\\
 Note that if $(X,d)$ is a geodesic space, then $d_\eps(x,y)=d(x,y)$ for all $\eps>0, x,y \in X$.
}\end{definition}

Throughout this paper, we consider a \emph{complete}, locally compact separable metric space $(X,d)$,
equipped with a Radon measure $m$ with full support, i.e., a Borel measure $m$ on $X$
which is finite on any compact set and strictly positive on any non-empty open set.
Such a triple $(X,d,m)$ is referred to as a \emph{metric measure space}. Then we set
$\diam(X,d):=\sup_{x,y\in X}d(x,y)$ and $B(x,r):=\{y\in X\mid d(x,y)<r\}$ for $x\in X$ and $r>0$.

Let $(\mathcal{E},\mathcal{F})$ be a \emph{symmetric Dirichlet form} on $L^{2}(X,m)$.
In other words, the domain $\mathcal{F}$ is a dense linear subspace of $L^{2}(X,m)$, such that
$\mathcal{E}:\mathcal{F}\times\mathcal{F}\to\mathbb{R}$
is a non-negative definite symmetric bilinear form which is \emph{closed}
($\mathcal{F}$ is a Hilbert space under the inner product $\sE_{1}(\cdot,\cdot):= \sE(\cdot,\cdot)+ \langle \cdot,\cdot \rangle_{L^{2}(X,m)}$)
and \emph{Markovian} (the unit contraction operates on $\sF$; $(u \vee 0)\wedge 1\in\mathcal{F}$ and $\mathcal{E}((u \vee 0)\wedge 1,(u \vee 0)\wedge 1)\leq \mathcal{E}(u,u)$ for any $u\in\mathcal{F}$).
Recall that $(\mathcal{E},\mathcal{F})$ is called \emph{regular} if
$\mathcal{F}\cap\mathcal{C}_{\mathrm{c}}(X)$ is dense both in $(\mathcal{F},\mathcal{E}_{1})$
and in $(\mathcal{C}_{\mathrm{c}}(X),\|\cdot\|_{\mathrm{sup}})$.
Here $\mathcal{C}_{\mathrm{c}}(X)$ is the space of $\mathbb{R}$-valued continuous functions on $X$
with compact support.

For a function $u \in \sF$, let $\supp_{m}[u]$ denote the support of the measure $|u|\,dm$,
i.e., the smallest closed subset $F$ of $X$ with $\int_{X\setminus F}|u|\,dm=0$; note that
$\supp_{m}[u]$ coincides with the closure of $X\setminus u^{-1}(\{0\})$ in $X$ if $u$ is continuous.
 Recall that
$(\mathcal{E},\mathcal{F})$ is called \emph{strongly local}
 if
$\mathcal{E}(u,v)=0$ for any $u,v\in\mathcal{F}$ with $\supp_{m}[u]$, $\supp_{m}[v]$
compact and $v$ is constant $m$-almost everywhere in a neighborhood of $\supp_m[u]$.
The pair $(X,d,m,\mathcal{E},\mathcal{F})$ of a metric measure space $(X,d,m)$ and a strongly local,
regular symmetric Dirichlet form $(\mathcal{E},\mathcal{F})$ on $L^{2}(X,m)$ is termed a \emph{metric measure Dirichlet space}, or an \textbf{MMD space}.
We refer to \cite{FOT,CF} for a comprehensive account of the theory of symmetric Dirichlet forms.

We recall the definition of energy measures associated to an MMD space. Note that $fg\in\mathcal{F}$
for any $f,g\in\mathcal{F}\cap L^{\infty}(X,m)$ by \cite[Theorem 1.4.2-(ii)]{FOT}
and that $\{(-n)\vee(f\wedge n)\}_{n=1}^{\infty}\subset\mathcal{F}$ and
$\lim_{n\to\infty}(-n)\vee(f\wedge n)=f$ in norm in $(\mathcal{F},\mathcal{E}_{1})$
by \cite[Theorem 1.4.2-(iii)]{FOT}.

\begin{definition}[cf.\ {\cite[(3.2.14) and (3.2.15)]{FOT}}]\label{d:EnergyMeas} {\rm
	Let $(X,d,m,\mathcal{E},\mathcal{F})$ be an MMD space.
	The \textbf{energy measure} $\Gamma(f,f)$ of $f\in\mathcal{F}$ is defined,
	first for $f\in\mathcal{F}\cap L^{\infty}(X,m)$ as the unique ($[0,\infty]$-valued)
	Borel measure on $X$ with the property that
	\begin{equation}\label{e:EnergyMeas}
	\int_{X} g \, d\Gamma(f,f)= \mathcal{E}(f,fg)-\frac{1}{2}\mathcal{E}(f^{2},g) \qquad \mbox{ for all $g \in \mathcal{F}\cap\mathcal{C}_{\mathrm{c}}(X)$,}
	\end{equation}
	and then by
	$\Gamma(f,f)(A):=\lim_{n\to\infty}\Gamma\bigl((-n)\vee(f\wedge n),(-n)\vee(f\wedge n)\bigr)(A)$
	for each Borel subset $A$ of $K$ for general $f\in\mathcal{F}$. 
	
	The notion of energy measure can be extended to the local Dirichlet space $\sF_{\loc}$, which is defined as
	\begin{equation}\label{e:Floc}
	\mathcal{F}_{\loc} := \biggl\{ f \in L^2_{\loc}(X,m) \biggm|
	\begin{minipage}{260pt}
	For any relatively compact open subset $V$ of $X$, there exists
	$f^{\#} \in \mathcal{F}$ such that $f \one_{V} = f^{\#} \one_{V}$ $m$-a.e.
	\end{minipage}
	\biggr\}.
	\end{equation}
	For any $f \in \sF_{\loc}$ and for any relatively compact open set $V \subset X$, we define
	\[
	\Gam(f,f)(V)= \Gam(f^{\#},f^{\#})(V),
	\]	
	where $f^{\#}$ is as in the definition of $\sF_{\loc}$. Since $(\sE,\sF)$ is strongly local, the value of $\Gam(f^{\#},f^{\#})(V)$ does not depend on the choice of $f^{\#}$, and is therefore well defined. Since $X$ is locally compact, this defines a Radon measure $\Gam(f,f)$ on $X$.
}\end{definition}

\begin{definition}[Capacity between sets]{\rm
 For subsets $A,B \subset X$, we define
 \[
 \sF(A,B):= \set{f \in \sF: f \equiv 1 \mbox{ on a neighborhood of $A$ and } f \equiv 0 \mbox{ on a neighborhood of $B$}},
 \]
 and the capacity $\Cap(A,B)$ as 
\[
\Cap(A,B)= \inf \set{\sE(f,f): f \in \sF(A,B)}.
\]
}\end{definition}

The main result of this work provides upper bounds on $d_\eps(x,y)$ based on analytic conditions on an MMD space that we introduce now. 
 Henceforth, we fix a continuous increasing bijection $\Psi:(0,\infty) \to (0,\infty)$  such that for all $0 < r \le R$,
 \be  \label{e:reg}
 C^{-1} \left( \frac R r \right)^{\beta_1} \le \frac{\Psi(R)}{\Psi(r)} \le C \left( \frac R r \right)^{\beta_2}, 
 \ee
 for some constants $0 < \beta_1 < \beta_2$ and $C>1$. Throughout this work, the function $\Psi$ is meant to denote the space time scaling of the process.
\begin{definition}
	{\rm
We recall the following properties that an MMD space $(X,d,m,\sE,\sF)$ may satisfy:\\
We say that $(X,d,m)$ satisfies the volume doubling property \hypertarget{vd}{$\on{(VD)}$} if there exists $C_D \ge 1$ such that
\be \tag*{($\on{VD}$)} 
m(B(x,2r)) \le C_D m(B(x,r)), \q \mbox{ for all $x \in X, r >0$.}
\ee
	 We say that $(X,d,m,\sE,\sF)$ satisfies the \emph{Poincar\'e inequality} \hypertarget{pi}{$\operatorname{PI}(\Psi)$}, if there exist constants $C,A  \ge 1$ such that 
	 for all $x\in X$, $r \in (0, \infty)$ and $f \in \sF$
	 \be  \tag*{$\operatorname{PI}(\Psi)$}
	 \int_{B(x,R)} (f - \ol f)^2 \,dm  \le C \Psi(r) \, \int_{B(x,A r)}d\Gamma(f,f),
 \ee
	 where $\ol f= m(B(x,r))^{-1} \int_{B(x,r)} f\, d\mu$.
	 If all balls are relatively compact, then the above inequality can be extended to all $f \in \sF_{\loc}$.

 We say that  $(X,d,m,\sE,\sF)$ satisfies the \emph{capacity estimate} \hypertarget{cap}{$\operatorname{cap}(\Psi)_\le$} 	
 if there exist $C_1,A_1,A_2>1$ such that for all $x \in X$, $0<R< \diam(X,d)/A_2$,
 \be \tag*{$\operatorname{cap}(\Psi)_\le$}
 \Cap(B(x,R),B(x,A_1R)^c) \le C_1 \frac{m(B(x,R))}{\Psi(R)}. % C_1^{-1} \frac{\mu(B(x,R)}{\Psi(R)} \le
 \ee
 }
\end{definition}
We recall the definition of the heat kernel corresponding to an MMD space.
\begin{definition}{\rm
	Let $(X,d,m,\mathcal{E},\mathcal{F})$ be an MMD space, and let $\set{P_t}_{t>0}$
	denote its associated Markov semigroup. A family $\set{p_t}_{t>0}$ of non-negative
	Borel measurable functions on $X \times X$ is called the
	\emph{heat kernel} of $(X,d,m,\mathcal{E},\mathcal{F})$, if $p_t$ is the integral kernel
	of the operator $P_t$ for any $t>0$, that is, for any $t > 0$ and for any $f \in L^2(X,m)$,
	\[
	P_{t} f(x) = \int_X p_{t}(x,y) f (y)\, dm (y) \qquad \mbox{for $m$-almost every $x \in X$.}
	\]
We remark that not every MMD space  $(X,d,m,\mathcal{E},\mathcal{F})$  has a heat kernel. The existence of heat kernel is an issue in general.
}\end{definition}
\subsection{Main results}
Our main result is the following upper bound on $d_\eps$.
\begin{theorem} \label{t:main}
	Let $(X,d,m,\sE,\sF)$ be an MMD space that satisfies  \hyperlink{vd}{$\on{(VD)}$}, \hyperlink{pi}{$\operatorname{PI}(\Psi)$}, and \hyperlink{cap}{$\operatorname{cap}(\Psi)_\le$}, where $\Psi$ satisfies \eqref{e:reg}.
	Then there exists $C>1$ such that for all $\eps>0$ and for all $x,y \in X$ that satisfy $d(x,y) \ge \eps$, we have
	\be \label{e:main}
\frac{d_\eps(x,y)^2}{\eps^2} \le C \frac{\Psi(d(x,y))}{	\Psi(\eps) }
	\ee
	In particular, for all $x,y \in X$, we have
	\be \label{e:GT}
	\lim_{\eps \to 0} \Psi(\eps) \frac{d_\eps(x,y)}{\eps} = 0.
	\ee 
\end{theorem}
\begin{remark} \label{r:main} {\rm
	\begin{enumerate}[(a)]
		\item  If $\Psi(r)=r^2$,  \eqref{e:main} implies the chain condition $d_\eps(x,y) \lesssim d(x,y)$ for all $\eps>0, x,y \in X$.
		\item If $\Psi(r) = r^\beta$, then \eqref{e:main} and  the triangle inequality $d_\eps(x,y) \ge d(x,y)$ imply that
		\[
		\frac{d(x,y)^2}{\eps^2} \le \frac{d_\eps(x,y)^2}{\eps^2}  \le C \frac{d(x,y)^\beta}{\eps^\beta},
		\]
		for all $x,y \in X, \eps>0$ with $d(x,y) \ge\eps$. By letting $\eps \to 
		0$, we give a new proof of the known fact that $\beta \ge 2$ must  necessarily hold.
		\item Let $\Psi(r)=r^\beta$ with $\beta \ge 2$. Consider the Dirichlet form corresponding to the Brownian motion on $\bR^d$ with Lebesgue measure as the symmetric measure, and the snowflake metric $d(x,y)=  \norm{x-y}_2^{2/\beta}$, where $\norm{x-y}_2$ denotes the Euclidean distance  (cf. \cite[Defintion 1.2.8]{MT} for the terminology `snowflake metric'). In this case, it is easy to obtain 
		\[
		\Psi(\eps) \frac{d_\eps(x,y)^2}{\eps^2} \asymp  \Psi(d(x,y)) \asymp \norm{x-y}_2^2
		\]
		for all $x,y \in \bR^d, \eps> 0$ with $\eps < d(x,y)$. Hence, the bound  \eqref{e:main} is sharp  for all $\beta \ge 2$.
		\item Theorem \ref{t:main} provides a new proof to an estimate due to Hambly and Kumagai \cite[Lemma 3.3]{HaKu}. Based on the results in  \cite{HaKu}, the estimate \eqref{e:GT} was introduced by Grigor'yan and Telcs in \cite[(1.8)]{GT12} to obtain sharp estimates of  the heat kernel (cf. Corollary \ref{c:hke}).
	\end{enumerate}
}\end{remark}
By \cite[Theorem 6.5]{GT12} along with \eqref{e:GT}, we have the following corollary (see Theorem \ref{t:chain} for generalization to arbitrary scale functions $\Psi$).
\begin{cor} \label{c:chain}
	Let $(X,d,m,\sE,\sF)$ be an MMD space that satisfies the following sub-Gaussian estimate on its heat kernel $p_t(\cdot,\cdot)$: there exists $\beta \ge 2,C,c>0$ such that
	\be  \label{e:sge}
 \frac{c}{V(x,t^{1/\beta})} \exp \left( - \left(\frac{d(x,y)^\beta}{ct}\right)^{1/(\beta-1)}\right) \le	p_t(x,y) \le \frac{C}{V(x,t^{1/\beta})} \exp \left( - \left(\frac{d(x,y)^\beta}{Ct}\right)^{1/(\beta-1)}\right)
	\ee
	for all $x,y \in X, t>0$, where $V(x,r)=m(B(x,r))$.
	Then  the metric $d$ satisfies the chain condition: there exists $K>1$ such that 
	\be  \label{e:chain}
	 d(x,y) \le d_\eps(x,y) \le K d(x,y) \q \mbox{for all $\eps>0$ and for all $x,y \in X$}.
	\ee 
\end{cor}
\begin{remark}{\rm 
		\begin{enumerate}[(a)]
	 \item The chain condition \eqref{e:chain} admits the following characterization. Let  $(X,d)$ be metric space such that the open balls $B(x,r)$ are relatively compact for all $x \in X,r >0$. Then $(X,d)$ satisfies the chain condition if and only if there exists a geodesic metric $\rho$ such that $d$ is bi-Lipschitz equivalent to $\rho$ \cite[Proposition A.1]{KM}. %; that is, there exists $C>1$ such that $C^{-1} d \le \rho \le C d$.
	 Recall that $(X,\rho)$ is \emph{geodesic}, if for any two points $x,y \in X$, there exists a function $\gam:[0,\rho(x,y)] \to X$ such that $\rho(\gam(s),\gam(t))=\abs{s-t}$ for all $s,t \in [0,\rho(x,y)]$.
  	\item Corollary \ref{c:chain} was previously known only for the case $\beta=2$. By a version of Varadhan's  asymptotic formula in \cite[Theorem 1.1]{HR}, we obtain that $d$ is bi-Lipschitz equivalent to the intrinsic metric. Hence by the remark above, $d$ satisfies the chain condition. However, the intrinsic metric vanishes identically for the case $\beta>2$. This suggests the need for a completely different approach when $\beta>2$.
  	\item The chain condition plays an essential role in the proof of singularity of energy measures in \cite{KM} for spaces satisfying the sub-Gaussian heat kernel estimate.
\end{enumerate}
}\end{remark}

Recall that the parameter $\beta$ in Corollary \ref{c:hke} is called the \emph{walk dimension}.
The following result can be viewed as a generalization of the result that the walk dimension is at least two.
\begin{cor} \label{c:ge2}
		Let $(X,d,m,\sE,\sF)$ be an MMD space that satisfies  \hyperlink{vd}{$\on{(VD)}$}, \hyperlink{pi}{$\operatorname{PI}(\Psi)$}, and \hyperlink{cap}{$\operatorname{cap}(\Psi)_\le$}, where $\Psi$ satisfies \eqref{e:reg}. Then there exists $C_1 \ge 1$ such that
		\[
		\frac{\Psi(r)}{\Psi(s)} \ge C_1^{-1} \left(\frac{r}{s}\right)^2, \qq \mbox{for all $0<s \le r < \diam(X,d)$.}
		\]
\end{cor}
%By Corollary \ref{c:ge2}, we can assume that $2 \le \beta_1 \le \beta_2$ in \eqref{e:reg}. 
{\bf Notation}. In the following, we will use the notation  $A \lesssim B$ for quantities $A$ and $B$ to indicate the existence of an
implicit constant $C \ge 1$ depending on some inessential parameters such that $A \le CB$. We write $A \asymp B$, if $A \lesssim B$ and $B \lesssim A$.

\section{Proofs}
\subsection{Connectedness}

We recall the notion of a net in a metric space. 
\begin{definition} \label{d:net} {\rm 
Let $(X,d)$ be a metric space and let $\eps>0$. A maximal $\eps$-separated subset $V \subset X$ is called an $\eps$-net; in other words, $V$ satisfies the following properties:
\begin{enumerate}
	\item[(a)] $V$ is $\eps$-separated; that is,  $d(x,y) \ge \eps$  whenever $x,y \in V$ and $x \neq y$.
	\item[(b)](maximality) If $W \supseteq V$ and $W$ is $\eps$-separated, then $W=V$.
\end{enumerate}
}\end{definition}
% Our first result is $\eps$-connectedness under volume doubling and Poincar\'e inequality.
As a first step towards \eqref{e:main}, we show that $d_\eps(x,y)$ is finite for all $\eps >0, x,y \in X$.
\begin{lem} \label{l:connect}
	Let $(X,d,m,\sE,\sF)$ be an MMD space that satisfies \hyperlink{vd}{$\on{(VD)}$}, \hyperlink{pi}{$\operatorname{PI}(\Psi)$}, where $\Psi$ satisfies \eqref{e:reg}. Then 
	\[
	d_\eps(x,y) < \infty \q \mbox{ for all $x,y \in X$  and for all $\eps >0$.}
	\]
\end{lem}
\proof
Let $x \in X, \eps >0$ and  $U_x= \set{y: d_\eps(x,y)< \infty}, V_x= X \setminus U_x$.
Clearly $U_x$ is open, since $z \in U_x$ implies $B(z,\eps) \subset U_x$. Similarly, $V_x$ is open since $z \in V_x$ implies $B(z,\eps) \subset V_x$. The above statements also imply that
\be \label{e:cn1}
\inf \set{d(y,z) : y \in U_x, z \in V_x}  \ge \eps.
\ee
Let $N$ denote a $\eps/2$-net in $(X,d)$.
Let $\phi_z \in C_c(X) \cap \sF$ such that $1 \ge \phi_z \ge 0, \restr{\phi_z}{B(z,\eps/2)} \equiv 1$, and $\supp(\phi_z) \subset B(z,\eps)$. 
Define
\[
\phi(y)= \sup_{z \in N \cap U_x} \phi_z(y) \q \mbox{ for all $y \in X$}.
\]
By \eqref{e:cn1}, and  $\cup_{z \in N} B(z,\eps/2) = X$, we obtain
\be \label{e:cn2}
\phi \equiv 1_{U_x}.
\ee
For any precompact open set $U$, 
we have $ \abs{\set{z \in N \cap U_x: B(z,\eps) \cap U } } < \infty$.
In this case, by setting $N_U= \set{z \in N \cap U_x : B(z,\eps) \cap U }$, we have
\[
\inf_{z \in N_U} \phi_z(y)= \phi(y) \q \mbox{ for all $y\in U$}, \q \inf_{z \in N_U} \phi_z \in \sF \cap C(X).
\]
Therefore $\phi=1_{U_x} \in \sF_{\on{loc}}\cap C(X)$.
By \cite[Theorem 4.3.8]{CF}, the push-forward measure of $\Gam(\phi,\phi)$ by $\phi$ 
is absolutely continuous with respect to the $1$-dimensional Lebesgue measure.
Since $\set{0,1}$ has zero Lebesgue measure, we obtain
\be  \label{e:cn3}
\Gamma(\phi,\phi)(B)= 0
\ee
for all balls $B=B(x,R), R>0$.
Since $\phi$ is continuous, by \hyperlink{pi}{$\on{PI}(\Psi)$}, $\phi$ is constant on all balls $B(x,R)$. Therefore $\phi(y)= \phi(x)=1$ for all $y \in B(x,R), R>0$. Therefore $\phi \equiv 1$. Hence $V_x = \emptyset$. Since  $x \in X, \eps >0$ are arbitrary, we obtain the desired estimate.
\qed

Recall that a metric space $(X,d)$ is said to be \emph{uniformly perfect}, if there exists $C >1$ such that $B(x,r) \setminus B(x,r/C) \neq \emptyset$ for all $x \in X, r>0$ that satisfy $B(x,r) \neq X$.

\begin{cor}
Let $(X,d,m,\sE,\sF)$ be an MMD space that \hyperlink{vd}{$\on{(VD)}$}, \hyperlink{pi}{$\operatorname{PI}(\Psi)$}, where $\Psi$ satisfies \eqref{e:reg}. Then $(X,d)$ is a uniformly perfect metric space.	
\end{cor}
\proof
Let $B(x,r)$ be a ball such that $B(x,r)\neq X$. Let $y \in X \setminus B(x,r)$. By Lemma \ref{l:connect}, there exists an $\eps$-chain  $\set{x_i}_{i=0}^N, x_0=x, x_N=y$ between $x$ and $y$ for $\eps=r/4$.
Since $d(x,x_0)=0, d(x,x_N) \ge r$ and $\abs{d(x,x_i)-d(x,x_{i+1})} \le d(x_i,x_{i+1})\le r/4$ for all $i=0,\ldots,N-1$, we have $d(x,x_{j}) \in [r/2,r)$ for some $1 \le j \le N-1$.
Hence \be \label{e:up}B(x,r) \setminus B(x,r/2) \neq \emptyset. \ee
\qed
\subsection{A two point estimate using Poincar\'e inequality}
For two measures $m, \nu$ on $(X,d)$, for $R> 0, x \in X$, we define a `truncated maximal function'
\be 
M_R^m \nu(x)= \sup_{0< r < R} \frac{\nu(B(x,r))}{m(B(x,r))}.
\ee
If $\nu \ll m$, then the above expression is the truncated maximal function of the Radon-Nikodym derivative $\frac{d\nu}{dm}$. However, in the lemma below $\nu$ will be the energy measure, and hence the measure $\nu$ and $m$ might be mutually singular. In the following lemma, by $\sC(X)$ we mean the space of continuous functions on $X$.
\begin{lem} \label{l:poin} (see \cite[Lemma 5.15]{HeKo})
	Let $(X,d,m,\sE,\sF)$ be an MMD space that satisfies \hyperlink{vd}{$\on{(VD)}$}, \hyperlink{pi}{$\operatorname{PI}(\Psi)$}, where $\Psi$ satisfies \eqref{e:reg}.
	There exists $C>0$ such that for all $x_0 \in X, R>0$, $x, y \in B(x_0, C^{-1} R)$, and for all $u \in \sC(X) \cap \sF_{\loc}$
	\[
	\abs{u(x)-u(y)}^2 \le C \Psi(R) \left(M_R^m \Gamma(u,u)(x)+ M_R^m \Gamma(u,u)(y)\right),
	\]
	where $\Gam(u,u)$ denotes the energy measure of $u$.
\end{lem} 
\proof
The proof in \cite{HeKo} applies to our setting with minor modifications. 
For the convenience of the reader, we recall the proof below.

Let $\delta \in (0,1)$ be a constant that will be chosen later. 
For a ball $B$, by $u_B$ we denote $\frac{1}{m(B)} \int u\,dm$.
Define $B_{x,i}= B(x,2^{-i} \delta R)$ for $i \in \set{0,1,2,\ldots}$. We estimate 
\begin{align}
\abs{u(x)-u_{B_{x,0}}} &\le \sum_{i=0}^\infty \abs{u_{B_{x,i}}-u_{B_{x,i+1}}} = \sum_{i=0}^{\infty} \frac{1}{m(B_{x,i+1})}  \abs{\int_{B_{x,i+1}} \left(u-u_{B_{x,i}}\right)\,dm} \nonumber \\ 
& \le\sum_{i=0}^{\infty} \frac{1}{m(B_{x,i+1})}\int_{B_{x,i+1}} \abs{u-u_{B_{x,i}}}\,dm   \nonumber \\
&\lesssim \frac{1}{m(B_{x,i})}\int_{B_{x,i}} \abs{u-u_{B_{x,i}}}\,dm \q \mbox{ (since $B_{x,i+1}\subset B_{x,i}$ and by  \hyperlink{vd}{$\on{(VD)}$})} \nonumber \\
& \le \sum_{i=0}^{\infty} \left( \frac{1}{m(B_{x,i})}\int_{B_{x,i}} \abs{u-u_{B_{x,i}}}^2\,dm  \right)^{1/2}  \q \q \mbox{(by H\"older inequality)}\nonumber \\
&\lesssim \sum_{i=0}^\infty \Psi(2^{-i}\delta R)^{1/2} \left( \frac{1}{m(B_{x,i})} \int_{ B(x, A \delta 2^{-i}R)}  d\Gam(u,u) \right)^{1/2} \mbox{(by \hyperlink{pi}{$\on{PI(\Psi)}$})}\nonumber \\
&\lesssim \sum_{i=0}^\infty \Psi(R)^{1/2} 2^{-\beta_1 i/2} \left(M_{A \delta R}^m \Gamma(u,u)(x) \right)^{1/2} \q \q \mbox{ (by \eqref{e:reg} and \hyperlink{vd}{$\on{(VD)}$})}\nonumber \\
&\lesssim   \left( \Psi(R)M_{A \delta R}^m \Gamma(u,u)(x) \right)^{1/2}. \label{e:pn-a}
\end{align}
Similarly, by setting $B_{y}=B(y,\delta' R)$, we obtain
\be  \label{e:pn-b}
\abs{u(x)-u_{B_y}} \lesssim \left( \Psi(R)M_{A \delta' R}^m \Gamma(u,u)(y)\right)^{1/2}.
\ee
We choose $C>1$ large enough  and $\delta' > \delta$ so that $x,y \in B(x_0,C^{-1}R)$ implies $B_x \subset B_y$. For example, it suffices to choose $C, \delta, \delta'$ so that
\be  \label{e:pn1}
2C^{-1}+\delta < \delta'.
\ee
We also require
\be  \label{e:pn2}
A \delta \le A \delta' \le 1,
\ee
so that $M_{A \delta R}^m \Gamma(u,u)(x) \le M_R^m \Gam(u,u)(x)$ and $M_{A \delta R}^m \Gamma(u,u)(y) \le M_R^m \Gam(u,u)(y)$. 

Evidently, given any $A>1$, it is possible to choose $C>1, 0 < \delta \le \delta'<1$ such that \eqref{e:pn1} and \eqref{e:pn2} are satisfied.
We fix one such $C, \delta, \delta'$ for the remainder of the proof.

By $B_{x,0} \subset B_y$,  $m(B_y) \lesssim m(B_{x,0})$, \hyperlink{vd}{$\on{(VD)}$} and \hyperlink{pi}{$\on{PI(\Psi)}$}, we obtain
\begin{align}\label{e:pn-c}
\abs{u_{B_{x,0}}- u_{B_y}} &\le \frac{1}{m(B_{x,0})} \int_{B_{x,0}} \abs{u-u_{B_y}} \,dm \lesssim \frac{1}{m(B_y)} \int_{B_y} \abs{u-u_{B_y}} \,dm \nonumber \\
 &\lesssim  \left( \Psi(R) \frac{\Gamma(u,u)( B(y,A \delta' R))}{m( B_y)}\right)^{1/2} \le \ \left( \Psi(R)M_{A \delta' R}^m \Gamma(u,u)(y)\right)^{1/2}
\end{align}
Combining \eqref{e:pn-a}, \eqref{e:pn-b}, \eqref{e:pn-c}, $A \delta \le A \delta' \le 1$, we obtain the desired bound
\begin{align*}
\abs{u(x)-u(y)} &\le \abs{u(x)-u_{B_{x,0}}} + \abs{u_{B_{x,0}}-u_{B_y}} + \abs{u(y)-u_{B_y}}
\\
&\lesssim \Psi(R)^{1/2}  \left(M_R^m \Gamma(u,u)(x)+ M_R^m \Gamma(u,u)(y)\right)^{1/2}.
\end{align*}
\qed

The telescoping sum argument has been applied in the context of anomalous diffusions previously in \cite[p. 1654]{BCK}. However, in \cite{BCK} the argument is used in the `strongly recurrent case'. One of the main novelties of this work is to extract useful estimates from that argument without further assumptions on volume growth.
\subsection{A partition of unity}
The use of partition of unity with functions of small Dirichlet energy indexed by a net is well-known \cite[p. 235]{Kan},  \cite[p. 504]{BBK}. Since we do not have a reference to give for the requirement in (c) below, we provide the details.
\begin{lem} \label{l:partition}
	Let $(X,d,m,\sE,\sF)$ be an MMD space that satisfies \hyperlink{vd}{$(\on{VD})$}, and \hyperlink{cap}{$\operatorname{Cap}(\Psi)_\le$}. Let $\eps>0$ and let $V$ denote any $\eps$-net. Let $\eps<\diam(X,d)/A_2$, where $A_2 \ge 1$ is the constant in  \hyperlink{cap}{$\operatorname{Cap}(\Psi)_\le$}.
	 Then, there exists a family of functions $\set{\psi_z: z\in V}$ that satisfies the following properties:
	\begin{enumerate}[(a)]
		\item  $\set{\psi_z: z\in V}$ is partition of unity $\sum_{z \in V} \psi_z  \equiv 1$.
		\item $\psi_z \in C_c(X) \cap \sF$ with  $0 \le \psi_z \le 1$, $\restr{\psi_z}{B(z,\eps/4)} \equiv 1$, and $\restr{\psi_z}{B(z,5\eps/4)^c} \equiv 0$.
		\item For all $z \in V$, $z' \in V \setminus \set{z}$, we have $\restr{\psi_{z'}}{B(z,\eps/4)} \equiv 0$. (this follows from (a) and (b)).
		\item There exists $C>1$ such that for all $z\in V$,
		\[
		\sE(\psi_z,\psi_z) \le C \frac{m(B(z,\eps))}{\Psi(\eps)}.
		\]
	\end{enumerate}
\end{lem}
\proof
For $z \in V$, we define the corresponding `Voronoi cell' $R_z$ as
\[
R_z =\set{p \in X: d(p,z)= d(p,V)= \min_{v \in V} d(p,v)}.
\]
We denote its $\eps/4$-neighbourhood by
$
R_z^{\eps/4}= \cup_{x \in R_z} B(x,\eps/4).
$
By the triangle inequality, we have
\be \label{e:net1}
B(z,\eps/2) \subset R_z \subset \overline B(z,\eps), \q  \bigcup_{z \in V} R_z = X,
\ee
and
\be \label{e:net2}
p \in B(z,\eps/2) \mbox{ implies that }  p \notin V_{w} \mbox{ for $w \in V \setminus \set{z}$.}
\ee
By  \eqref{e:net2}, and the triangle inequality, we have
\be \label{e:net3}
v,w \in V \mbox{ and } v \neq w, \mbox{ imply that } B(z,\eps/4) \cap R_w^{\eps/4} = \emptyset.
\ee
For $z \in V$, let $N_z$ denote an $\eps/(4A_1)$-net of $R_z$, where $A_1$ denote the constant in \hyperlink{cap}{$\on{cap}(\Psi)_\le$}. 
%By the volume doubling property and \eqref{e:net1}, 
%\be
%\sup_{z \in V} \abs{N_z}  < \infty.
%\ee
For each $w \in N_z$, by \hyperlink{cap}{$\on{cap}(\Psi)_\le$}, \eqref{e:reg}, \hyperlink{vd}{$\on{(VD)}$}, there exists $C_1>1$ such that for $w \in N_z, z \in V$, we have a non-negative function $\rho_w \in C_c(X) \cap \sF$ that satisfies
\be \label{e:net4}
\restr{\rho_w}{B(w,\eps/(4A_1)} \equiv 1,  \q \restr{\rho_w}{B(w,\eps/4)^c} \equiv 0, \q \sE(\rho_w,\rho_w) \le C_1 \frac{m(B(w,\eps)}{\Psi(\eps)}.
\ee
Hence by \hyperlink{vd}{$\on{(VD)}$} and \eqref{e:net4}, we obtain a  family of functions $\set{\phi_z: z \in V}$ that satisfy
\be  \label{e:net5}
\phi_z = \max_{w \in N_z} \rho_w \mbox{ such that } \restr{\phi_w}{R_z} \equiv 1, \restr{\phi_w}{(R_z^{\eps/4})^c} \equiv 0, \sE(\phi_z,\phi_z) \lesssim \frac{m(B(z,\eps)}{\Psi(\eps)}.
\ee
Define 
\[
\psi_z:= \frac{\phi_z}{ \sum_{w \in V} \phi_w}.
\]
Property (a) is immediate. Properties (b) and (c) follow from \eqref{e:net5}, \eqref{e:net3}, and \eqref{e:net1}.
Property (d) follows from  \eqref{e:net5}, \hyperlink{vd}{$\on{(VD)}$},  Leibniz rule, chain rule and Cauchy-Schwarz inequality, we obtain
\begin{align*}
\sE(\psi_z,\psi_z) &\lesssim \sup \left(  \sum_{w \in V}  \phi_w \right)^{-2} \left( \sE(\phi_z,\phi_z) + (\sup \phi_z)^2 \sum_{w \in V \cap B(z,5 \eps/2)}  \sE(\phi_w,\phi_w) \right)\\
&\lesssim \sum_{w \in V \cap B(z,5\eps/2)} \frac{m(B(w,\eps))}{\Psi(\eps)} \lesssim  \frac{m(B(z,\eps))}{\Psi(\eps)}.
\end{align*}
\qed
\subsection{Proof of the main result}
We recall an elementary lemma from \cite{GT12}. 
\begin{lemma}  (\cite[Lemma 6.3]{GT12}) \label{l:chain}
	Let $(X,d)$ be a metric space.
 Define $N_\eps(x,y)$ as the minimal value of $N$ such that there exists an $\eps$-chain $\set{x_i}_{i=0}^N$ between $x$ and $y$. If $d_\eps(x,y) < \infty$ for some $x,y \in X, \eps>0$, then 
 \[
 \left\lceil \frac{d_\eps(x,y)}{\eps} \right\rceil \le N_\eps(x,y) \le 9 \left\lceil \frac{d_\eps(x,y)}{\eps} \right\rceil.
 \]	
\end{lemma}

\noindent{\em Proof of Theorem \ref{t:main}.}
%\subsection{An application to heat kernel estimates}
Let $A_2$ denote the constant in  \hyperlink{cap}{$\operatorname{Cap}(\Psi)_\le$}.
Since $d_{\eps} \le d_{\eps'}$ whenever $\eps' \le \eps$, by replacing $\eps$ by $\eps/(2A_2)$ if necessary and by using \eqref{e:reg}, we assume that $\eps< \diam(X,d)/A_2$.

Fix $x,y \in X, \eps>0$ such that $d(x,y) \ge \eps$. Set $\eps'=\eps/3$. Let $V$ be an $\eps'$-net such that $\set{x,y} \subset V$. 
Define $\hat{u}: V \to [0,\infty)$ as 
\be \label{e:mn1}
\hat{u}(z):= N_\eps(x,z),
\ee
where $N_\eps(x,z)$ is as defined in Lemma \ref{l:chain}. 
By Lemma \ref{l:connect}, $\hat{u}$ is finite. By definition,
\be \label{e:mn2}
\abs{\hat{u}(z_1)- \hat{u}(z_2)} \le 1 \q \mbox{ for all $z_1,z_2 \in V$ such that $d(z_1,z_2)<\eps$.}
\ee
Let   $\set{\psi_z: z\in V}$  denote the partition of unity defined in Lemma \ref{l:partition}.
Define  $u:X \to [0,\infty)$ as
\bes
u(p):= \sum_{z \in V} \hat u(z)\psi_z(p). 
\ees
For any ball $B(x_0,r), x_0 \in X, r >0$, by Lemma \ref{l:partition} we have
\be \label{e:mn3}
u(p) = \sum_{z \in V \cap B(x_0, r+ 5\eps'/4)}  \hat u(z)\psi_z(p)\q \mbox{ for all $p \in B(x_0,r)$.}
\ee
Since $V \cap B(x_0, r+ 5\eps'/4)$ is a finite set by \hyperlink{vd}{$\on{(VD)}$}, we obtain that $u \in \sF_{\loc}$.
By Lemma \ref{l:partition}(b), we have $\restr{u}{B(z,\eps'/4)}\equiv \hat u(z)$ for all $z \in V$. Therefore by \cite[Theorem 4.3.8]{CF}, 
the push-forward measure of $\Gam(u,u)$ by $u$ is absolutely continuous with respect to the $1$-dimensional Lebesgue measure.
Therefore, we obtain
\be \label{e:mn4}
\Gamma(u,u)(B(z,\eps'/4)) =0 \q \mbox{ for all $z \in V$.} 
\ee
By \eqref{e:mn3} and Lemma \ref{l:partition}(a), we have
\be \label{e:mn5}
u(p)= \hat u (z)+ \sum_{w \in V \cap B(z, 9\eps'/4)}  (\hat u(w)- \hat u(z))\psi_w(p)\q \mbox{ for all $p \in B(z,\eps'), z \in V$.}
\ee
By \hyperlink{vd}{$\on{(VD)}$}, there exits $C_1>1$ such that $\sup_{z \in V} \abs{V \cap B(z, 9\eps'/4)} \le C_1$. By \eqref{e:mn5}, and Cauchy-Schwarz inequality, there exists $C_2>1$ such that the following holds: for all $z \in V$, we have
\begin{align} \label{e:mn6}
\Gamma(u,u)(B(z,\eps')) &\le C_1 \sum_{w \in V \cap B(z, 9\eps'/4)} (\hat u(w)- \hat u(z))^2 \sE(\psi_w,\psi_w) \nonumber \\
&\lesssim \sum_{w \in V \cap B(z, 9\eps'/4)}  \frac{m(B(w,\eps'))}{\Psi(\eps')} \q \mbox{ (by \eqref{e:mn2} and Lemma \ref{l:partition}(d))} \nonumber \\
&\le C_2 \frac{m(B(z,\eps'/2))} {\Psi(\eps)} \q \mbox{ (by \hyperlink{vd}{$\on{(VD)}$} and \eqref{e:reg})}.
\end{align}
By Lemma \ref{l:partition}, \eqref{e:mn6} and \hyperlink{vd}{$\on{(VD)}$}, there exists $C_3>0$ such that for all $z \in X, r \ge \eps'/4$, we have
\be \label{e:mn7}
\Gamma(u,u)(B(z,r)) \le  C_2 \sum_{w \in B(z,r+5\eps'/4)}  \frac{m(B(w,\eps'/2))} {\Psi(\eps)} \le C_3 \frac{m(B(z,r))}{\Psi(\eps)}.
\ee
Combining \eqref{e:mn4} and \eqref{e:mn7}, we obtain
\be \label{e:mn8}
M^m_R\Gamma(u,u)(z)= \sup_{r <R} \frac{\Gamma(u,u)(B(z,r))}{m(B(z,r))}\le  \frac{C_3}{\Psi(\eps)} \q \mbox{for all $z \in V, R>0$.}
\ee 
By \eqref{e:mn8}, Lemma \ref{l:poin}, $\hat u (x)=0, \hat u(y)= N_\eps(x,y)$, and \eqref{e:reg}, there exists $C_4 >0$ such that
\bes
N_\eps(x,y)^2 \le C_4 \frac{\Psi(d(x,y))}{\Psi(\eps)} \mbox{ for all $x,y \in X , \eps \le d(x,y)$.}
\ees
Combining the above estimate along with Lemma \ref{l:chain}, we obtain \eqref{e:main}.  We obtain \eqref{e:GT} using \eqref{e:main} and $\lim_{\eps \downarrow 0} \Psi(\eps)=0$.
\qed
\begin{remark}{\rm
\begin{enumerate}[(a)]
	\item The constant $C$ in \eqref{e:main} can be chosen to depend only on the constants associated with the assumptions \hyperlink{vd}{$\on{(VD)}$}, \hyperlink{pi}{$\operatorname{PI}(\Psi)$},  \hyperlink{cap}{$\operatorname{cap}(\Psi)_\le$} and \eqref{e:reg}.
	\item The proof of the estimate on $d_\eps(x,y)$ uses \hyperlink{cap}{$\operatorname{cap}(\Psi)_\le$} only at scales less than $\eps$, whereas it relies on \hyperlink{pi}{$\operatorname{PI}(\Psi)$} for scales up to the order of $d(x,y)$. In other words, our argument relies on the Poincar\'e inequality on a larger range of scales than it relies on the capacity upper bound. 
\end{enumerate}	
}\end{remark}

\noindent{\em Proof of Corollary \ref{c:ge2}.}
Let  $0<s < r < \diam(X,d)$. By  uniform perfectness of $(X,d)$ (more precisely, by \eqref{e:up}), there exists $x,x_r \in X$ such that $r/2 \le d(x,x_r)<r$.

By \eqref{e:reg}, it suffices to consider the case $s \le r/2$. Therefore,
\begin{align*}
\frac{\Psi(r)}{\Psi(s)} &\asymp 
\frac{\Psi(d(x,x_r))}{\Psi(s)}  \qq \mbox{(by \eqref{e:reg})}\\
 & \gtrsim \left(\frac{d(x,x_r)}{s} \right)^2 \gtrsim \left(\frac r s\right)^2 \qq \mbox{(by \eqref{e:main}).}
\end{align*}

\subsection{Application to heat kernel estimates}
\begin{definition}{\rm
	%	Let $\Psi$ satisfy \eqref{e:reg}.
		
		Let $U \subset V$ be open sets in $X$ with $U \subset \overline{U} \subset V$. We say a continuous
		function $\phi$ is a cutoff function for $U \subset V$ if $\phi=1$ on a neighbourhood of $\overline{U}$ and $\supp(\phi) \subset V$.\\
		We recall the cutoff Sobolev inequality  \hypertarget{cs}{$\operatorname{CS}(\Psi)$} for an MMD space$(X,d,m,\sE,\sF)$: there exists $C_S>0$ such that for any $x \in X, R,r>0$, there exists a cutoff function $\phi$ for $B(x,R) \subset B(x,R+r)$ such that
		\[
		\int_{B(x,R+r) \setminus B(x,R)} f^2\, d\Gam(\phi,\phi) \le \frac{1}{8} \int_{B(x,R+r) \setminus B(x,R)} \phi^2 \, d\Gam(f,f) + \frac{C_S}{\Psi(r)} \int_{B(x,R+r) \setminus B(x,R)} f^2\,dm,
		\]
		for all $f \in \sF$.}
\end{definition}

For $\Psi$ satisfying \eqref{e:reg}, we define
\be \label{e:defPhi}
\Phi(s)= \sup_{r>0} \left(\frac{s}{r}-\frac{1}{\Psi(r)}\right).
\ee

\begin{cor} \label{c:hke}
	Let $(X,d,m,\sE,\sF)$ be an MMD space that satisfies  \hyperlink{vd}{$\on{(VD)}$}, \hyperlink{pi}{$\operatorname{PI}(\Psi)$}, and \hyperlink{cs}{$\operatorname{CS}(\Psi)$}, where $\Psi$ satisfies \eqref{e:reg}.
	Then we have matching upper and lower bounds on the heat kernel as given in 
	\be \label{e:hke}
	p_t(x,y) \asymp \frac{C}{V(x,\Psi^{-1}(t))} \exp \left(-c t \Phi\left(\frac{d_\eps(x,y)}{t} \right)\right),
	\ee
	where $\Phi$ is as defined in \eqref{e:defPhi}, and $\eps=\eps(t,x,y)$ is chosen so that
	\[
	\eps(t,x,y)= \sup \set{\eps>0 : \frac{\Psi(\eps)}{\eps}d_\eps(x,y) \le t }.
	\]
	Here $\asymp$ in \eqref{e:hke} means that both $\le$ and $\ge$ are true, but the positive constants $C$ and $c$ may be different for upper and lower bounds.
\end{cor}
\proof 
This follows from \cite[Theorem 1.2]{GHL}, \cite[Theorem 6.5]{GT12} and Theorem \ref{t:main}.
\qed
\begin{lem} \label{l:legendre}
	Let $\Psi:(0,\infty) \to (0,\infty)$ be a continuous, increasing bijection such that 
	\be \label{e:Psireg}
	C_1^{-1} \left( \frac R r \right)^{\beta_1} \le \frac{\Psi(R)}{\Psi(r)} \le C_1 \left( \frac R r \right)^{\beta_2}, 
	\ee
	for some constants $1 < \beta_1 \le \beta_2$ and $C_1>1$.
	Let $\Phi$ be the function defined in 
	\eqref{e:defPhi}.
	Then, there exists $C_2>0$ such that 
	\be \label{e:Phireg}
	C_2^{-1} \left( \frac S s \right)^{\beta_2/(\beta_2-1)} \le \frac{\Phi(S)}{\Phi(s)} \le C_2 \left( \frac S s \right)^{\beta_1/(\beta_1-1)}, \q \mbox{for all $0 < s \le S$.}
	\ee
\end{lem}
\proof
For $\lambda \ge 1, r>0$, we have
\begin{align*}
\Phi(s) &=\sup_{r>0} \left(\frac{s}{r}-\frac{1}{\Psi(r)}\right) = \sup_{r>0} \left(\frac{s}{r \lam}-\frac{1}{\Psi(r \lam)}\right) \\
& \le \sup_{r>0} \left(\frac{s}{r \lam}-\frac{1}{C_1 \lam^{\beta_2}\Psi(r)}\right)\q \mbox{(by \eqref{e:Psireg})}\\
&\le \frac{1}{C_1 \lam^{\beta_2}}  \sup_{r>0} \left(\frac{s C_1^{-1}\lam^{\beta_2-1}}{r}-\frac{1}{\Psi(r)}\right) \\
& \le  \frac{1}{C_1 \lam^{\beta_2}} \Phi(s C_1^{-1}\lam^{\beta_2-1}).
\end{align*}
By choosing $S/s=C_1^{-1}\lam^{\beta_2-1}$ in the above estimate, we obtain the lower bound on $\frac{\Phi(S)}{\Phi(s)}$ in \eqref{e:Phireg}. The proof of upper bound in \eqref{e:Phireg} is analogous.
\qed

We prove the following generalization of Corollary \ref{c:chain}.  
\begin{thm} \label{t:chain}
	Let $(X,d,m,\sE,\sF)$ be an MMD space. Let $\Psi:(0,\infty)\to (0,\infty)$ be an increasing bijection that satisfies  \eqref{e:Psireg}.  Assume that $(X,d,m,\sE,\sF)$  satisfies the following estimate on
	 its heat kernel $p_t(\cdot,\cdot)$:
	\be  \label{e:hksge}
	p_t(x,y) \asymp \frac{C}{V(x,\Psi^{-1}(t))} \exp \left(-c t \Phi\left(\frac{d(x,y)}{t} \right)\right),
	\ee
	for all $x,y \in X, t>0$, where $V(x,r)=m(B(x,r))$, where  $\Phi:(0,\infty)\to (0,\infty)$ be as defined in \eqref{e:defPhi}.
		Here $\asymp$ in \eqref{e:hksge} means that both $\le$ and $\ge$ are true, but the positive constants $C$ and $c$ may be different for upper and lower bounds.
	Then  the metric $d$ satisfies the chain condition: there exists $K>1$ such that 
	\bes
	d(x,y) \le d_\eps(x,y) \le K d(x,y) \q \mbox{for all $\eps>0$ and for all $x,y \in X$}.
	\ees
\end{thm}
\proof
By \cite[Lemma 6.4]{GT12} and \eqref{e:Psireg}, $\eps=\eps(t,x,y)$ satisfies
\be \label{e:ch1}
\eps^{\beta_2 -1} d_\eps(x,y) \lesssim \frac{\Psi(\eps)}{\eps}d_\eps(x,y) = t.
\ee
Note that 
\begin{align}
t\Phi\left(\frac{2\Psi^{-1}(t)}{t} \right) &= \sup_{r>0} \left(\frac{2 \Psi^{-1}(t)}{ r} - \frac{t}{\Psi(r)}\right)  \nonumber \\
& \le \sup_{r>0} \left(\frac{2 \Psi^{-1}(t)}{ r} - C_1^{-1}\min\left(\frac{\Psi^{-1}(t)^{\beta_1}}{r^{\beta_1}},\frac{\Psi^{-1}(t)^{\beta_2}}{r^{\beta_2}}\right)\right) \q \mbox{(by \eqref{e:Psireg})} \nonumber \\
& \lesssim 1 \label{e:ch2}
\end{align}
 For $x \in X, r>0$, by integrating the lower bound in \eqref{e:sge}, over the ball $B(x,2r)$ with $\Psi^{-1}(t)=r$ and using \eqref{e:ch2}, we obtain
\[
\frac{V(x,2r)}{V(x,r)} \lesssim \int_{B(x,2r)} p_t(x,y)\,m(dy) \le 1,
\]
which implies the volume doubling property.
By the argument in the proof of Theorem 1.2 in \cite[proof of  (8.6)]{GHL}, the MMD space $(X,d,m,\sE,\sF)$ satisfies  \hyperlink{pi}{$\operatorname{PI}(\Psi)$}.
By \eqref{e:up} and \cite[Exercise 13.1]{Hei}, $(X,d,m)$ satisfies the reverse volume doubling property; that is, there exist $c,\alpha>0$ such that 
\[
\frac{V(x,R)}{V(x,r) }\ge c^{-1} \left(\frac R r\right)^\alpha,\ \q \mbox{for all $0<r<R <\diam(X,d)$.}
\]
By \cite[Theorem 1.2]{GHL}, the MMD space $(X,d,m,\sE,\sF)$ satisfies 
\hyperlink{vd}{$\on{(VD)}$}, \hyperlink{pi}{$\operatorname{PI}(\Psi)$}, and \hyperlink{cs}{$\operatorname{CS}(\Psi)$}. In \cite[Theorem 1.2]{GHL} assumes that the space is unbounded, the same result also extends to the bounded case \cite[Theorem 3.2]{Lie}; see \cite[Remark 2.9]{KM} for further discussion.
Therefore by Corollary \ref{c:hke} and \eqref{e:hksge}, for all $x,y \in X$ and $t>0$, we have
\begin{align*}
c_2\exp\left(-C_2 t \Phi\left(\frac{d(x,y)}{t} \right)\right) &\le V(x,\Psi^{-1}(t)) p_t(x,y) \\
&\le {C_1} \exp \left(-c_1 t \Phi\left(\frac{d_\eps(x,y)}{t} \right)\right)
\end{align*}
where $\eps=\eps(t,x,y)$.
Therefore for $x,y \in X$, $t>0$
\bes 
c_1 t \Phi\left(\frac{d_\eps(x,y)}{t}\right)  \le   C_2 t \Phi\left(\frac{d(x,y)}{t}\right) + \log(C_1/c_2)
\ees
Hence, by \cite[Lemma 3.19]{GT12}, for every $x,y \in X$, there exists $t_{x,y}>0$ such that for all $0<t<t_{x,y}$, we have
\bes
c_1 t \Phi\left(\frac{d_\eps(x,y)}{t}\right)  \le   2C_2 t \Phi\left(\frac{d(x,y)}{t}\right) 
\ees
By the above estimate and Lemma \ref{l:legendre}
\be \label{e:cc1}
d_\eps(x,y) \lesssim d(x,y), \q\mbox{for all $x,y >0, t< t_{x,y}$, where $\eps=\eps(t,x,y)$.}
\ee
By \eqref{e:ch1} and $d_\eps \ge d$, we obtain 
\[
\eps(t,x,y) \lesssim  \left(\frac{t}{d(x,y)} \right)^{1/(\beta_2-1)}.
\]
Hence  $\lim_{t\to 0}\eps(t,x,y)=0$.
Since $\eps \mapsto d_\eps(x,y)$ is non-increasing, and $\lim_{t\downarrow 0} \eps(t,x,y)=0$, by letting $t \downarrow 0$ in \eqref{e:cc1}, we obtain
\[
\lim_{\eps \downarrow 0} d_\eps(x,y) \lesssim d(x,y), \q \mbox{for all $x,y \in X$.}
\]

%Hence $d_{\eps(t,x,y)}(x,y) \asymp d(x,y)$ for all $x,y \in X, t>0$.
Hence for all $x,y \in X, \eps >0$, we obtain $d_\eps(x,y) \asymp d(x,y)$. 
\qed

\noindent{\em Proof of Corollary \ref{c:chain}.}
This is a special case of Theorem \ref{t:chain} with $\Psi(r)=r^\beta$.
  \qed

\begin{remark}{\rm
Let $\Psi:X \times [0,\infty) \to [0,\infty)$ is a regular scale function in the sense of \cite[Definition 5.4]{BM}; that is, there exists $C>1, \beta_1, \beta_2 >0$ such that for all $x,y \in X, 0 < s \le r$, we have
\[
C^{-1} \left(\frac{r}{d(x,y) \vee r}\right)^{\beta_2}\left(\frac{d(x,y) \vee r}{s}\right)^{\beta_1} \le  \frac{\Psi(x,r)}{\Psi(y,s)} \le C \left(\frac{r}{d(x,y) \vee r}\right)^{\beta_1}\left(\frac{d(x,y) \vee r}{s}\right)^{\beta_2}.
\]
Further, assume that the MMD space $(X,d,m,\sE,\sF)$ satisfies \hyperlink{vd}{$\on{(VD)}$}, \hyperlink{pi}{$\operatorname{PI}(\Psi)$}, and \hyperlink{cs}{$\operatorname{CS}(\Psi)$}.
In this case,
we can  obtain sharp heat kernel bounds as follows: by a change of metric as done in \cite[Proposition 5.7]{BM}, we can reduce it to the case $\Psi(r)=r^\beta$, which can be handled using Corollary \ref{c:hke}.}
\end{remark}
\subsection*{Acknowledgements} 
I thank Takashi Kumagai for inspiring this work by asking if Corollary \ref{c:chain} could  be true.
It is a pleasure to acknowledge  Naotaka Kajino and Takashi Kumagai for helpful references and  discussions.
This work was carried out at the Research
Institute for Mathematical Sciences, Kyoto University and at Kobe University. I am grateful to Takashi Kumagai,  Naotaka Kajino, Jun Kigami,  and
Ryoki Fukushima  for support and hospitality during the visit. 
The author thanks the anonymous referee for a careful reading of the paper and helpful suggestions.

\noindent Department of Mathematics, University of British Columbia,
Vancouver, BC V6T 1Z2, Canada. \\
mathav@math.ubc.ca

\end{document}